\documentclass{arx}
\usepackage{amsmath,amssymb,url}
\usepackage{enumitem} 
\usepackage[all]{xy} 
\usepackage{xcolor}
\theoremstyle{plain}
\newtheorem{theorem}{Theorem}[section]
\newtheorem{lemma}[theorem]{Lemma}

\newtheorem{corollary}[theorem]{Corollary}

\theoremstyle{definition}
\newtheorem{definition}[theorem]{Definition}

\numberwithin{equation}{section}

\theoremstyle{definition}

\newtheorem{problem}[theorem]{Problem}

\numberwithin{equation}{section}

\begin{document}

\title[On the Boolean algebra tensor product]{On the Boolean algebra tensor product via Carath\'{e}odory spaces of place functions}

\author{Gerard Buskes}
\address{Department of Mathematics, University of Mississippi, University, MS 38677}
\email{mmbuskes@olemiss.edu}

\author{Page Thorn}
\address{Department of Mathematics, University of Mississippi, University, MS 38677}
\email{lthorn@go.olemiss.edu}

\subjclass{46A40 , 46M05, 06E99}

\date{January 27, 2023}


\keywords{Riesz space, Vector lattice, Boolean algebra, Tensor product, Free product, Dedekind complete}

\begin{abstract} We show that the Carath\'{e}odory space of place functions on the free product of two Boolean algebras is Riesz isomorphic with Fremlin's Archimedean Riesz space tensor product of their respective Carath\'{e}odory spaces of place functions. We provide a solution to Fremlin's problem 315Y(f) in  \cite{fremlin_measure} concerning completeness in the free product of Boolean algebras by applying our results on the Archimedean Riesz space tensor product to Carath\'{e}odory spaces of place functions. 
\end{abstract}

\maketitle

\section{Introduction and preliminary material}

Fremlin asserts in problem 315Y(f) of \cite{fremlin_measure} that the Boolean algebra tensor product of two nontrivial Boolean algebras is complete if and only if one is finite and the other is complete. In Theorem 4.6 of \cite{buskes_thorn_complete}, we prove that the Fremlin tensor product of two Dedekind complete Riesz spaces rarely is Dedekind complete. In fact, if the tensor product is Dedekind complete, then one of the two spaces is Riesz isomorphic to the set of all finite-valued functions on a subset of that space. To connect 315Y(f) of \cite{fremlin_measure} with Theorem 4.6 of \cite{buskes_thorn_complete}, we employ Carath\'{e}odory's Riesz space of place functions on a Boolean algebra. The main result is Theorem \ref{CA main result} with applications given in Section \ref{applications}.

The necessary terms for Boolean algebras, the free product, Riesz spaces, and Carath\'{e}odory spaces of place functions are provided in this section. We reserve $\mathcal{A}$, $\mathcal{B}$ for Boolean algebras and $E$, $F$, $G$ for Archimedean Riesz spaces.

\subsection*{Boolean algebras and their free product}

For Boolean algebras, see chapter 31 of \cite{fremlin_measure}. Two elements $x$ and $y$ of a Boolean algebra are called \emph{disjoint} if $x\wedge y =0$, in which case we write $x\perp y$. Two subsets $A$ and $B$ of a Boolean algebra are called disjoint if $x\perp y$ for every $x\in A$ and $y\in B$, in which case we write $A\perp B$. We define the \emph{disjoint sum} of two elements $x$ and $y$ in a Boolean algebra by 
$$x\oplus y=(x\wedge y')\vee (x' \wedge y).$$
A Boolean algebra is \emph{complete} if every nonempty subset has a supremum.

\begin{definition}(312F of \cite{fremlin_measure}) Let $\mathcal{A}$ and $\mathcal{B}$ be Boolean algebras. A function ${\chi\colon\mathcal{A}\to\mathcal{B}}$ is said to be a \emph{Boolean homomorphism} if for all $x$, $y\in\mathcal{A}$,
\begin{enumerate}[label=(\roman*)]
\item $\chi (x\wedge y)=\chi(x)\wedge\chi(y)$;
\item $\chi(x\oplus y)=\chi(x)\oplus\chi(y)$;
\item $\chi(1_\mathcal{A})=1_\mathcal{B}$. 
\end{enumerate}
A bijective Boolean homomorphism is called a \emph{Boolean isomorphism}. If there exists an isomorphism $\chi\colon\mathcal{A}\to\mathcal{B}$, then the Boolean algebras $\mathcal{A}$ and $\mathcal{B}$ are said to be \emph{Boolean isomorphic}.
\end{definition}

Proposition 312H of \cite{fremlin_measure} proves additionally that every Boolean homomorphism preserves finite suprema, that is $\chi(x\vee y)=\chi(x)\vee\chi(y)$ for every $x$, $y\in\mathcal{A}$. 

The \emph{Stone space} of a Boolean algebra $\mathcal{A}$ is the set $Z$ of nonzero ring homomorphisms from $\mathcal{A}$ to $\mathbb{Z}_2$. Set $$\hat{a}=\{z:z\in Z, z(a)=1\}.$$ By Stone's Theorem (see, for instance, 311E of \cite{fremlin_measure}), the canonical map $$a\mapsto\hat{a}:\mathcal{A}\to\mathcal{P}(\mathbb{Z})$$ is an injective ring homomorphism which we call the \emph{Stone representation}. For more on Stone spaces, see \cite{fremlin_measure} where Fremlin defines and utilizes Stone spaces to define the Boolean algebra tensor product, called the \emph{free product}.


\begin{definition}\label{free product}(Fremlin, 315I of \cite{fremlin_measure}) \begin{enumerate}[label=(\roman*)]
\item Let $\{\mathcal{A}_i\}_{i\in I}$ be a family of Boolean algebras. For each $i\in I$, let $Z_i$ be the Stone space of $\mathcal{A}_i$. Set $Z=\prod_{i\in I} Z_i$, with the product topology. Then the \emph{free product} of $\{\mathcal{A}_i\}_{i\in I}$ is the algebra of open-and-closed sets in $Z$, denoted $\bigotimes_{i\in I} \mathcal{A}_i$.
\item For $i\in I$ and $a\in\mathcal{A}_i$, the set $\hat{a}\subseteq Z_i$ representing $a$ is an open-and-closed subset of $Z_i$; because $z\mapsto z(i)\colon Z\to Z_i$ is continuous, $$\epsilon_i(a)=\{z:z(i)\in\hat{a}\}$$ is open-and-closed, so belongs to $\mathcal{A}$. In this context, $\epsilon_i\colon\mathcal{A}_i\to\mathcal{A}$ is called the \emph{canonical map}.
\end{enumerate}\end{definition}

In the following theorem, we list the necessary material from 315J and 315K of \cite{fremlin_measure} in the language of Fremlin.
\begin{theorem}\label{map properties} Let $\{\mathcal{A}_i\}_{i\in I}$ be a family of Boolean algebras, with free product $\mathcal{A}$. 
\begin{enumerate}[label=(\roman*)]
\item The canonical map $\epsilon_i\colon\mathcal{A}_i\to\mathcal{A}$ is a Boolean homomorphism for every $i\in I$.
\item For any Boolean algebra $\mathcal{B}$ and any family $\{\varphi_i\}_{i\in I}$ such that $\varphi_i$ is a Boolean homomorphism from $\mathcal{A}_i$ to $\mathcal{B}$ for every $i$, there is a unique Boolean homomorphism $\varphi\colon\mathcal{A}\to\mathcal{B}$ such that $\varphi_i=\varphi\circ\epsilon_i$ for each $i$.
\item Write $C$ for the set of those members of $\mathcal{A}$ expressible in the form $\inf_{j\in J}\epsilon_j(a_j)$, where $J\subseteq I$ is finite and $a_j\in\mathcal{A}_j$ for every $j$. Then every member of $\mathcal{A}$ is expressible as the supremum of a disjoint finite subset of $C$.
\item $\mathcal{A}=\{0_{\mathcal{A}}\}$ if and only if there is some $i\in I$ such that $\mathcal{A}_i=\{0_{\mathcal{A}_i}\}$.
\item If $\mathcal{A}_i\neq\{0_{\mathcal{A}}\}$ for every $i\in I$, then $\epsilon_i$ is injective for every $i\in I$.
\item Let $\mathcal{A}_i\neq\{0_{\mathcal{A}}\}$ for every $i\in I$. If $J\subseteq I$ is finite and $a_j$ is a nonzero member of $\mathcal{A}_j$ for each $j\in J$, then $\inf_{j\in J}\epsilon_j(a_j)\neq 0$.
\end{enumerate}
\end{theorem} 

\subsection*{Archimedean Riesz spaces and their tensor product}
See \cite{zaanen_introduction_1997} for Archimedean Riesz spaces and \cite{fremlin_tensor_1972} for Riesz bimorphisms.

\begin{theorem}\label{tensor product} (4.2 of \cite{fremlin_tensor_1972}) Let $E$ and $F$ be Archimedean Riesz spaces. There exists an Archimedean Riesz space $G$ and a Riesz bimorphism $\varphi\colon E\times F\to G$ such that whenever $H$ is an Archimedean Riesz space and $\psi\colon  E\times F\to H$ is a Riesz bimorphism, there is a unique Riesz homomorphism $T\colon G\to H$ such that $T\circ \varphi = \psi$.
\end{theorem}

$G$ of Theorem \ref{tensor product}  is the \emph{Archimedean Riesz space tensor product of $E$ and $F$}, denoted by $E\bar{\otimes}F$. The ``universal property of $E\bar{\otimes}F$" refers to the implication that any Archimedean Riesz space paired with a Riesz bimorphism satisfying Theorem \ref{tensor product} is Riesz isomorphic to $G$. The Riesz bimorphism $\otimes\colon E\times F \to E\bar{\otimes} F$ embeds the algebraic tensor product $E\otimes F$ into $E\bar{\otimes}F$ via $\otimes (e,f) = e\otimes f$ for all $e\in E$ and $f\in F$. 

We define a few terms needed for the statement of Theorem \ref{DC main result}.

\begin{definition} Let $E$ be an Archimedean Riesz space and let $I$ be a nonempty set. $c_{00}(I,E)$ is the set of all maps $f\colon I\to E$ for which $$S(f)=\{x\in I : f(x)\neq 0\}$$ is finite. If $E=\mathbb{R}$, then $c_{00}(I,E)$ is written $c_{00}(I)$. \end{definition}

Let $f$ and $g$ be elements of a Riesz space $E$. The ideals generated by $f$ and $g$ are denoted by $E_f$ and $E_g$ respectively. We denote the principal bands generated by $f$ and $g$ with $[f]$ and $[g]$ respectively (see, for instance, pg. 30 of \cite{zaanen_introduction_1997} for definitions).

A Riesz space is \emph{Dedekind complete} if every bounded subset of $E$ has a supremum. Every Dedekind complete Riesz space is Archimedean. In \cite{buskes_thorn_complete}, we characterized when the tensor product of two Dedekind complete Riesz spaces is Dedekind complete.

\begin{theorem}{\label{DC main result}}(4.6 of \cite{buskes_thorn_complete}) Suppose E and F are Dedekind complete Riesz spaces. The following are equivalent.
\begin{enumerate}
\item $E_x \bar{\otimes} F_y$ is Dedekind complete for every $x \in E^+$ and $y \in F^+$.
\item $[E_x$ is finite dimensional for every $x \in E^+]$ or [$F_y$ is finite dimensional for every $y \in F^+]$.
\item $E \cong c_{00}(I)$ for a set $I \subseteq E$ or $F \cong c_{00}(J)$ for a set $J \subseteq F$.
\item $E\bar{\otimes}F \cong c_{00}(I,F)$ for a set $I \subseteq E$ or $E\bar{\otimes}F \cong c_{00}(J,E)$ for a set $J \subseteq F$.
\item $E\bar{\otimes}F$ is Dedekind complete.
\end{enumerate}
\end{theorem}

As an intermediary between an Archimedean Riesz spaces and Boolean algebras, we consider Boolean algebras of bands. The following three statements are used in Section \ref{applications} and are given for the reader's convenience.

\begin{theorem}\label{complete zaanen}(22.6, 22.8 of \cite{luxemburg_riesz_1971})
Let $E$ be a Riesz space and define $$\mathcal{B}(E)=\{B\subseteq E : B\text{ is a band}\}.$$ $\mathcal{B}(E)$ is an order complete distributive lattice. $\mathcal{B}(E)$, partially ordered by inclusion, is a Boolean algebra if and only if $E$ is Archimedean. \end{theorem}

\begin{lemma}\label{disjoint bands} Let $E$ be a Riesz space and $f$, $g\in E$. Then $|f|\wedge |g| =0$ implies $[f]\perp [g].$  \end{lemma}
\begin{proof}
Let $|f|\wedge|g|=0$. Certainly, $E_f\perp E_g$. Suppose $h_1\in [f]$ and $h_2\in [g]$. Using the fact that $E_f\perp E_g$, it is straightforward to show that $|h_1|\wedge |h_2|=0$ via 6.1 and 7.8 of \cite{zaanen_introduction_1997}. 
Thus, $[g_1]\perp [g_2]$.
\end{proof}

\begin{lemma}\label{inf dim BA of bands} If $E$ is an infinite dimensional Archimedean Riesz space, then $\mathcal{B}(E)$ is not finite.\end{lemma}
\begin{proof}
By the contrapositive of Theorem 26.10 in \cite{luxemburg_riesz_1971}, there is an infinite subset of mutually disjoint nonzero elements in $E$. Thus, there are an infinite number of mutually disjoint bands in $E$ by Lemma \ref{disjoint bands}.
\end{proof}

\subsection*{Carath\'{e}odory spaces of place functions}

\begin{definition}(pg. 40 of \cite{aliprantis_positive_2006}) Let $E$ be a Riesz space and $e\in E ^+$. Then $x\in E^+$ is said to be a \emph{component} of $e$ whenever $x\wedge(e-x)=0$. \end{definition}
The collection of all components of $e$, denoted $\mathcal{C}(e)$, is a Boolean algebra under the partial ordering induced by $E$ (pg. 40 of \cite{aliprantis_positive_2006}). With $e$ as a strong order unit (pg. 51 of \cite{zaanen_introduction_1997}), a connection between Archimedean Riesz spaces and Boolean algebras is described explicitly in the following theorem. 
\begin{theorem}\label{definition of place functions}(4.1 of \cite{buskes_loomis_sikorski}) Let $\mathcal{A}$ be a Boolean algebra. There exists an Archimedean Riesz space $E$ with a strong unit $e$ with the following properties. \begin{enumerate}[label=(\roman*)]
\item There exists a Boolean isomorphism $\chi\colon\mathcal{A}\to\mathcal{C}(e)$.
\item $E$ is the linear span of $\mathcal{C}(e)$.
\end{enumerate}
\end{theorem}
$(E,\chi)$ is unique up to isomorphism. It is denoted by $\mathcal{C}(\mathcal{A})$ and is called the \emph{space of place functions on $\mathcal{A}$} in the sense of Carath\'{e}odory. 

Let $\lambda_i$, $\gamma_j \in \mathbb{R}$ be nonzero; $n$, $m \in \mathbb{N}$; $x_i \in \mathcal{A}$ be pairwise disjoint; and $y_j \in \mathcal{A}$ be pairwise disjoint. 
Two elements $$f= \sum_{i=1}^n \lambda_i \chi(x_i) \hspace{.5cm} \text{ and } \hspace{.5cm} g= \sum_{j=1}^m \gamma_j \chi(y_j)$$ are equivalent if $\bigvee_{i=1}^n x_i = \bigvee_{j=1}^m y_j$ and if $\lambda_i =\gamma_j$ whenever $x_i\wedge y_j\neq 0$. $\mathcal{C}(\mathcal{A})$ is the set of all such equivalence classes. Henceforth, we take $f= \sum_{i=1}^n \lambda_i \chi(x_i)$ to represent all elements of $\mathcal{C}(\mathcal{A})$ that are equivalent to $f$.
\par
We define addition in $C(\mathcal{A})$ in the style of Goffman in \cite{goffman_1958} and Jakubik in \cite{jakubik_place}. 
For a different approach, see \cite{buskes_loomis_sikorski}.
For $x$, $y \in \mathcal{A}$, let ${x}-_1 y$ be the complement of $x \wedge y$ relative to $x$, that is, $x\wedge (x \wedge y)'$. Then addition in $C(\mathcal{A})$ is defined by 
$$f+g
=\sum_{i=1}^n \sum_{j=1}^m (\lambda_i + \gamma_j) \chi(x_i\wedge y_j) 
+ \sum_{i=1}^n \lambda_i \chi(x_i-_1 \bigvee_{j=1}^m y_j)
+ \sum_{j=1}^m \gamma_j \chi(y_j-_1 \bigvee_{i=1}^n x_i)$$
where in the summation only those terms are taken into account in which $\lambda_i + \gamma_j \neq 0$ and the elements $x_i \wedge y_j$, $x_i-_1 \bigvee y_j$, and $y_j-_1 \bigvee x_i$ are non-zero. It is routine to verify that addition is well-defined in $\mathcal{C}(\mathcal{A})$.

Jakubik proves in \cite{jakubik_place} that the completeness of a Boolean algebra is equivalent to the Dedekind completeness of its Carath\'{e}odory space of place functions. However, his propositions assume complete distributivity. Since this work has no need for a Boolean algebra to be completely distributive, Theorem \ref{BA and CS complete} is proven with credit to Propositions 5.2(a) and 5.6 of \cite{jakubik_place} for its similarity.

\begin{definition}(page 231 of \cite{jakubik_place}) Let $Y$ be a sublattice of a lattice $X$. $Y$ is said to be a \emph{regular sublattice} of $X$ if:
\begin{enumerate}[label=(\roman*)]
\item whenever $x_0\in Y$ and $\emptyset\neq X\subseteq Y$ such that $x_0=\sup_{Y} X$, then $x_0=\sup_{X} X$; and 
\item whenever $x_1\in Y$ and $\emptyset\neq X\subseteq Y$ such that $x_1=\inf_{Y} X$, then $x_1=\inf_{X} X$.
\end{enumerate}
\end{definition}

\begin{theorem}\label{BA and CS complete} Let $\mathcal{A}$ be a Boolean algebra. $\mathcal{A}$ is complete if and only if $\mathcal{C}(\mathcal{A})$ is Dedekind complete. \end{theorem}
\begin{proof}
Assume that $\mathcal{A}$ is complete. Let $D$ be a bounded subset of $\mathcal{C}(\mathcal{A})$. Then there exists $g\in\mathcal{C}(\mathcal{A})$ such that $g\geq f$ for every $f\in\mathcal{C}(\mathcal{A})$.
Find $\lambda_i\in\mathbb{R}$, $n\in\mathbb{N}$, and $x_i\in\mathcal{A}$ such that $g = \sum_{i=1}^{n} \lambda_i \chi(x_i)$. Set $$x=x_1\vee\cdots\vee x_n \hspace{1cm} \text{ and } \hspace{1cm} \lambda = max\{\lambda_1,\cdots,\lambda_n\}.$$
Then $D\subseteq[0,\lambda \chi(x)]$.
By assumption, the interval $[0,x]$ is complete in $\mathcal{A}$.
It follows from Corollary 4.4 of \cite{jakubik_place} that  $\mathcal{A}$ is a regular subset of $\mathcal{C}(\mathcal{A})$. Then the interval $[0,\chi(x)]$ is complete as a subset of $\mathcal{C}(\mathcal{A})$.
In particular, $[0,\lambda \chi(x)]$ is complete, so $\sup (D)$ exists in $\mathcal{C}(\mathcal{A})$.
\par 
To prove sufficiency, assume that $\mathcal{C}(\mathcal{A})$ is complete. Let $\chi \colon \mathcal{A} \to \mathcal{C}(e)$ be the Boolean isomorphism from Theorem \ref{definition of place functions}. Note that $e=\chi(1_\mathcal{A})$.
Let $D$ be a subset of $\mathcal{A}$. 
Since $\mathcal{C}(\mathcal{A})$ is Dedekind complete, $\sup \chi (D)$ exists in $\mathcal{C}(\mathcal{A})^+$. For every $x\in D$, $\chi(x)$ is a component of $\chi(1_\mathcal{A})$. Thus, $\sup\chi (D) = 2\sup\chi (D) \wedge \chi(1_\mathcal{A})$ so that $0 = \sup\chi (D) \wedge \left(\chi(1_\mathcal{A})-\sup\chi (D)\right)$.
By definition, $\sup \chi(D)$ is a component of $e$.
\par 
Let $y=\chi^{-1}(\sup\chi (D))$.
Since $\chi$ is a Boolean isomorphism, $y$ is an upper bound for $D$.
Suppose there exists $y'$ such that $x\leq y' < y$ for every $x \in D$.
Then $\chi(y')\geq \sup_{x\in D} \chi (x) = \chi (y)$.
Thus, $\chi(y')=\chi(y)$. Since $\chi$ is one-to-one, $y'=y$. 
Therefore, $y=\sup (D)$ exists in $\mathcal{A}$.
\end{proof}

\section{The Fremlin tensor product of Carath\'{e}odory spaces of place functions}
In this section, we relate Boolean algebras $\mathcal{A}$, $\mathcal{B}$, and $\mathcal{A}\otimes\mathcal{B}$ to their Carath\'{e}odory spaces of place functions $\mathcal{C}(\mathcal{A})$, $\mathcal{C}(\mathcal{B})$, and $\mathcal{C}(\mathcal{A}\otimes\mathcal{B})$. The notation of Theorem \ref{definition of place functions} is used with the addition of subscripts to indicate which Boolean algebra is at work. The symbols in $(1)$, $(2)$, and $(3)$ will be used freely.
\begin{enumerate}[label=(\arabic*)]
\item $\chi_A\colon \mathcal{A}\to C(\mathcal{A})$, $\chi_B\colon \mathcal{B}\to C(\mathcal{B})$, and $\hat{\chi}\colon \mathcal{A} {\otimes} \mathcal{B} \to C(\mathcal{A} {\otimes} \mathcal{B})$ are the Boolean isomorphisms from Theorem \ref{definition of place functions}.
\item $C(\mathcal{A})$, $C(\mathcal{B})$, and $C(\mathcal{A} {\otimes} \mathcal{B})$ have units $\chi_A(1_\mathcal{A})$, $\chi_B(1_\mathcal{B})$, and $\hat{\chi}(1_{\mathcal{A}\otimes\mathcal{B}})$ respectively.
\item $\epsilon_A\colon\mathcal{A}\to\mathcal{A}\otimes\mathcal{B}$ and $\epsilon_B\colon\mathcal{B}\to\mathcal{A}\otimes\mathcal{B}$ are the canonical Boolean homomorphisms in Definition \ref{free product}.
\end{enumerate} 

\begin{theorem}\label{CA main result} $C(\mathcal{A}) \bar{\otimes} C(\mathcal{B})$ and $C(\mathcal{A} \otimes \mathcal{B})$ are Riesz isomorphic.\end{theorem}
\begin{proof}
Assume that $\mathcal{A}$ and $\mathcal{B}$ are nontrivial Boolean algebras.
For $f\in \mathcal{C}(\mathcal{A})$, there exist $n\in\mathbb{N}$, pairwise disjoint $x_i \in \mathcal{A}$, and nonzero $\lambda_i\in\mathbb{R}$ such that ${f= \sum_{i=1}^n \lambda_i \chi_A(x_i)}$. For $g\in  \mathcal{C}(\mathcal{B})$, there exist $m\in\mathbb{N}$, pairwise disjoint $u_j\in\mathcal{B}$, and nonzero $\gamma_j\in\mathbb{R}$ so $g= \sum_{j=1}^m \lambda_j \chi_B(u_j)$.
Define $\psi\colon \mathcal{C}(\mathcal{A})\times \mathcal{C}(\mathcal{B})\to \mathcal{C}(\mathcal{A} \otimes \mathcal{B})$ by
\begin{align*}
\psi(f,g)=& \psi\left(\sum_{i=1}^n \lambda_i \chi_A(x_i), \sum_{j=1}^m \gamma_j \chi_B(u_j)\right)\\
=& \sum_{i=1}^n \sum_{j=1}^m (\lambda_i \gamma_j) \hat{\chi}(\epsilon_A(x_i)\wedge \epsilon_B(u_j)).
\end{align*} 
It follows from Theorem \ref{map properties} (iv) and (vi) that the definition of $\psi$ is independent of the representations chosen for $f$ and $g$.
\par
Let $f_1=f$ and $f_2 = \sum_{k=1}^p \delta_k \chi_A(y_k)$ for nonzero $\delta_k\in \mathbb{R}$, $p\in \mathbb{N}$ and pairwise disjoint $y_k \in \mathcal{A}$. Recall that $f_1+f_2$ is defined to be $$\sum_i \sum_k (\lambda_i+\delta_k) \chi_A(x_i\wedge y_k) 
+ \sum_i \lambda_i \chi_A({x_i}_{-1} \bigvee_k y_k)
+ \sum_k \delta_k \chi_A({y_k}_{-1} \bigvee_i x_i).$$

\noindent \textbf{Claim}: $\psi$ is bilinear.
\begin{align*}\small
& \psi(f_1 + f_2,g)\\
=& \psi\left(\sum_i \lambda_i \chi_A(x_i) + \sum_k \delta_k \chi_A(y_k), \sum_j \gamma_j \chi_B(u_j)\right)\\
=&\sum_{i,k,j} (\lambda_i+\delta_k)\gamma_j\hat{\chi}\left(\epsilon_A(x_i\wedge y_k)\wedge \epsilon_B(u_j)\right) + \sum_{i,j} (\lambda_i\gamma_j) \hat{\chi}\left(\epsilon_A({x_i}_{-1} \bigvee_k y_k)\wedge \epsilon_B(u_j)\right)\\
& + \sum_{k,j} (\delta_k \gamma_j)\hat{\chi}\left(\epsilon_A({y_k}_{-1} \bigvee_i x_i)\wedge \epsilon_B(u_j)\right)\\
=&\sum_{i,j}\left[ \sum_k (\lambda_i\gamma_j)\hat{\chi}\left(\epsilon_A(x_i\wedge y_k)\wedge \epsilon_B(u_j)\right)+ (\lambda_i\gamma_j) \hat{\chi}\left(\epsilon_A({x_i}_{-1} \bigvee_k y_k)\wedge \epsilon_B(u_j)\right)\right]\\
& + \sum_{k,j} \left[ \sum_i(\delta_k \gamma_j)\hat{\chi}\left(\epsilon_A(x_i\wedge y_k)\wedge \epsilon_B(u_j)\right) + (\delta_k \gamma_j)\hat{\chi}\left(\epsilon_A({y_k}_{-1} \bigvee_i x_i)\wedge \epsilon_B(u_j)\right)\right]\\
=&\sum_{i,j}(\lambda_i\gamma_j)\left[\hat{\chi}\left(\bigvee_k\epsilon_A(x_i\wedge y_k)\wedge \epsilon_B(u_j)\right)+ \hat{\chi}\left(\epsilon_A({x_i}_{-1} \bigvee_k y_k)\wedge \epsilon_B(u_j)\right)\right]\\
& + \sum_{k,j}(\delta_k \gamma_j)\left[\hat{\chi}\left(\bigvee_i\epsilon_A(x_i\wedge y_k)\wedge \epsilon_B(u_j)\right) + \hat{\chi}\left(\epsilon_A({y_k}_{-1} \bigvee_i x_i)\wedge \epsilon_B(u_j)\right)\right]\hspace{.5cm}(*)\\
=&\sum_{i,j} (\lambda_i\gamma_j)\hat{\chi}(\epsilon_A(x_i)\wedge \epsilon_B(u_j))+ \sum_{k,j}(\delta_k \gamma_j)\hat{\chi}(\epsilon_A(y_k)\wedge \epsilon_B(u_j))\\
=& \psi(f_1,g) + \psi(f_2,g).
\end{align*}
\normalsize
$(*)$ is justified because $y_k\perp y_{k'}$ for all $k\neq k'$ and $x_i\perp x_{i'}$ for all $i\neq i'$.  Symmetrically, $\psi(f,g_1+g_2)=\psi(f,g_1)+\psi(f,g_2)$ for $f\in \mathcal{C}(\mathcal{A})$ and $g_1$, $g_2 \in \mathcal{C}(\mathcal{B})$. It follows from the definition of $\psi$ that $\psi(\lambda f,g)=\psi(f,\lambda g)= \lambda \psi(f,g)$ for every $\lambda\in\mathbb{R}$.

\noindent\textbf{Claim}: $\psi$ is a Riesz bimorphism.\\ 
Assume $f_1\wedge f_2=0$ and $g\in\mathcal{C}(\mathcal{B})^+$. Using the same representations for $f_1$, $f_2$, and $g$ as above, it follows that $x_i\perp y_k$ for all $i$ and $k$.
Then since the maps $\hat{\chi}$ and $\epsilon_A$ are Boolean homomorphisms and $\{x_i\}_{i=1}^n$, $\{y_k\}_{k=1}^p$ are each pairwise disjoint,
\begin{align*}
\psi(f_1,g)\wedge\psi(f_2,g)=& 
\psi\left(\underset{i}\sum \lambda_i \chi_A(x_i), \sum_j \gamma_j \chi_B(u_j)\right) \wedge \psi\left(\sum_k \delta_k \chi_A(y_k), \sum_j \gamma_j \chi_B(u_j)\right)\\
=& \left(\sum_{i,j}(\lambda_i \gamma_j) \hat{\chi}(\epsilon_A(x_i)\wedge \epsilon_B(u_j))\right) \wedge \left(\sum_{k,j} (\delta_k \gamma_j) \hat{\chi}(\epsilon_A(y_k)\wedge \epsilon_B(u_j))\right)\\
=& \ 0.
\end{align*} 
Likewise if $f\in\mathcal{C}(\mathcal{A})^+$ and $g_1\wedge g_2=0$ in $\mathcal{C}(\mathcal{B})$, then $\psi(f,g_1)\wedge\psi(f, g_2)=0$. By Theorem 19.1 of  \cite{zaanen_introduction_1997}, $\psi$ is a Riesz bimorphism.

It follows from the universal property of the Fremlin tensor product that there is a unique Riesz homomorphism $T\colon C(\mathcal{A}) \bar{\otimes} C(\mathcal{B})\to C(\mathcal{A} \otimes \mathcal{B})$ such that $\psi=T\circ\otimes$.
$$\xymatrix{
C(\mathcal{A}) \times C(\mathcal{B}) \ar[d]_{\psi} 
\ar[r]^{\otimes} & C(\mathcal{A}) \bar{\otimes} C(\mathcal{B})\ar@{.>}[ld]^{T}\\
C(\mathcal{A}\otimes \mathcal{B}) 
}$$
\noindent\textbf{Claim}: $T$ is a Riesz isomorphism.\\
\noindent\textbf{Step 1}: $T$ is onto.\\
Let $h\in C(\mathcal{A}\otimes \mathcal{B})$. Then $h=\sum_{i=1}^n \lambda_i \hat{\chi}(e_i)$ for some pairwise disjoint $e_i\in \mathcal{A}\otimes\mathcal{B}$, $n\in\mathbb{N}$, and nonzero $\lambda_i\in\mathbb{R}$. Fix $i\in\{1,\cdots,n\}$. By Theorem \ref{map properties} $(iii)$, there exists a finite disjoint subset $\{\epsilon_A(a_{k})\wedge\epsilon_B(b_k)\}_{k=1}^m$
($m\in\mathbb{N}$) of $\mathcal{A}\otimes\mathcal{B}$ such that $$e_i=\bigvee_{k=1}^m\epsilon_A(a_{k})\wedge\epsilon_B(b_k).$$
Then it follows from the definition of $\psi$ that 
\begin{align*}
\hat{\chi}(e_i)=& \hat{\chi}\left(\bigvee_{k=1}^m\epsilon_A(a_{k})\wedge\epsilon_B(b_k)\right)\\
=&\bigvee_{k=1}^m\hat{\chi}(\epsilon_A(a_{k})\wedge\epsilon_B(b_k))\\
=&\bigvee_{k=1}^m\psi(\chi_A(a_k),\chi_B(b_k))\\
=&\bigvee_{k=1}^mT\circ\otimes(\chi_A(a_k),\chi_B(b_k)).
\end{align*}
Since $T$ preserves finite suprema, $\hat{\chi}(e_i)$ is in the image of $T$ for every $i$. It follows from the linearity of $T$ that $h$ is in the image of $T$.

\noindent\textbf{Step 2}: $T$ is one-to-one.\\
Suppose $f\in C(\mathcal{A}) \otimes C(\mathcal{B})$, the algebraic tensor product of $C(\mathcal{A})$ and $C(\mathcal{B})$, such that $f$ is nonzero. Then for some $n\in\mathbb{N}$, nonzero $\lambda_k\in\mathbb{R}$, and nontrivial $x_k\in\mathcal{A}$, $u_k\in\mathcal{B}$ such that $$f=\sum_{k=1}^n \lambda_k \chi_A(x_k)\otimes \chi_B(u_k).$$
Since $\epsilon_A$, $\epsilon_B$, and $\hat{\chi}$ are injective Boolean isomorphisms,
\begin{align*}
T(f)=& T\left(\sum_{k=1}^n \lambda_k \chi_A(x_k)\otimes \chi_B(u_k)\right)\\
=& \sum_{k=1}^n \lambda_k \psi\left(\chi_A(x_k),\chi_B(u_k)\right)\\
=& \sum_{k=1}^n \lambda_k \hat{\chi}(\epsilon_A(x_k)\wedge \epsilon_B(u_j))\\
\neq& 0.
\end{align*}
Let $g$ be a nonzero element of the Riesz space tensor product of $C(\mathcal{A})$ and $C(\mathcal{B})$, i.e. $g\in C(\mathcal{A}) \bar{\otimes} C(\mathcal{B})$. By Theorem 2.2 of \cite{allenby_labuschagne}, for all $\delta > 0$ there exists $f\in\mathcal{C}(\mathcal{A})^+\otimes \mathcal{C}(\mathcal{B})^+$ such that $0\leq |g|-f\leq \delta \hat{\chi}(1_{\mathcal{A}\otimes\mathcal{B}})$. Since $C(\mathcal{A}) \bar{\otimes} C(\mathcal{B})$ is Archimedean, choose $\delta >0$ such that $|g|\wedge \delta  \hat{\chi}(1_{\mathcal{A}\otimes\mathcal{B}})\neq |g|$. Then $f$ is nonzero. We have shown that $T(f)\neq 0$ when $0\neq f\in \mathcal{C}(\mathcal{A})\otimes \mathcal{C}(\mathcal{B})$. Since $T$ is a Riesz homomorphism, $0< T(f)\leq |T(g)|$. Therefore, $T(g)\neq 0$, and $T$ is a Riesz isomorphism. Consequently, $C(\mathcal{A}) \bar{\otimes} C(\mathcal{B})$ is Riesz isomorphic to $C(\mathcal{A} \otimes \mathcal{B})$.
\end{proof}

\section{Applications}\label{applications}

In this section, we use Theorem \ref{CA main result} to provide a solution for Fremlin's problem 315Y(f) in \cite{fremlin_measure}. The statement leads to an observation on Dedekind completeness in the Fremlin tensor product of place functions and a statement on bands in the Fremlin tensor product of infinite dimensional Archimedean Riesz spaces.

\begin{problem}\label{Fremlin complete BA}(Fremlin, 315Y(f) of  \cite{fremlin_measure}) Let $\mathcal{A}$ and $\mathcal{B}$ be Boolean algebras. $\mathcal{A}\otimes\mathcal{B}$ is complete if and only if \emph{either} $\mathcal{A}=\{0\}$ \emph{or} $\mathcal{B}=\{0\}$ \emph{or} $\mathcal{A}$ is finite and $\mathcal{B}$ is complete \emph{or} $\mathcal{B}$ is finite and $\mathcal{A}$ is complete. \end{problem}
\begin{proof} If $\mathcal{A}=\{0\}$ or $\mathcal{B}=\{0\}$, the result is trivial. Assume $\mathcal{A}$ and $\mathcal{B}$ are nontrivial Boolean algebras.

Suppose $\mathcal{A}\otimes\mathcal{B}$ is complete.
It follows from Theorems \ref{CA main result} and \ref{BA and CS complete} that $\mathcal{C}(\mathcal{A}\otimes\mathcal{B})\cong\mathcal{C}(\mathcal{A})\bar{\otimes}\mathcal{C}(\mathcal{B})$  is Dedekind complete. 
By Proposition 3.6 of \cite{grobler_2022}, $\mathcal{C}(\mathcal{A})$ and $\mathcal{C}(\mathcal{B})$ are Dedekind complete. From Theorem \ref{BA and CS complete}, $\mathcal{A}$ and $\mathcal{B}$ are complete. It remains to show that one of the Boolean algebras is finite. However, the Dedekind completeness of $\mathcal{C}(\mathcal{A})\bar{\otimes}\mathcal{C}(\mathcal{B})$ implies that $\mathcal{C}(\mathcal{A}) \cong c_{00}(I)$ for a set $I \subseteq \mathcal{C}(\mathcal{A})$ or $\mathcal{C}(\mathcal{B}) \cong c_{00}(J)$ for a set $J \subseteq \mathcal{C}(\mathcal{B})$ (see Theorem \ref{DC main result}). Since each Carath\'{e}odory space of place functions contains a unit, $\mathcal{C}(\mathcal{A})$ or $\mathcal{C}(\mathcal{B})$ is finite dimensional. Thus, $\mathcal{A}$ is finite or $\mathcal{B}$ is finite.
\par
The sufficiency is proven analogously via Theorem \ref{DC main result}.
\end{proof}

\begin{corollary} Let $\mathcal{A}$ and $\mathcal{B}$ be nontrivial Boolean algebras. $\mathcal{C}(\mathcal{A})\bar{\otimes}\mathcal{C}(\mathcal{B})$ is Dedekind complete if and only if one of $\mathcal{A}$ or $\mathcal{B}$ is finite and the other is complete. 
\end{corollary}


Recall that for an Archimedean Riesz space $E$, its collection of bands, denoted $\mathcal{B}(E)$, forms a complete Boolean algebra. Our last application shows that for  Archimedean Riesz spaces $E$ and $F$, the set of bands in $E\bar{\otimes}F$ is rarely Boolean isomorphic to $\mathcal{B}(E)\otimes\mathcal{B}(F)$. That is, if $E$ and $F$ are infinite dimensional, not every band $B$ of $E\bar{\otimes}F$ can be ``decomposed" into the Fremlin tensor product of a band in $E$ and a band in $F$.

\begin{corollary} Let $E$ and $F$ be infinite dimensional Archimedean Riesz spaces. Then $\mathcal{B}(E)\otimes\mathcal{B}(F)$ is not Boolean isomorphic to $\mathcal{B}(E\bar{\otimes} F)$.\end{corollary}
\begin{proof} By Lemma \ref{inf dim BA of bands}, neither $\mathcal{B}(E)$ nor $\mathcal{B}(F)$ is finite.
Then $\mathcal{B}(E)\otimes\mathcal{B}(F)$ is not complete by Problem \ref{Fremlin complete BA}. However, Theorem \ref{complete zaanen} states that the Boolean algebra of bands is complete for any Archimedean Riesz space, so $\mathcal{B}(E\bar{\otimes} F)$ is complete.
\end{proof}

\bibliographystyle{amsplain}
\bibliography{myreferences}

\newcommand{\noop}[1]{}
\providecommand{\bysame}{\leavevmode\hbox to3em{\hrulefill}\thinspace}
\providecommand{\MR}{\relax\ifhmode\unskip\space\fi MR }
\providecommand{\MRhref}[2]{%
  \href{http://www.ams.org/mathscinet-getitem?mr=#1}{#2}
}
\providecommand{\href}[2]{#2}
\begin{thebibliography}{10}

\bibitem{aliprantis_positive_2006}
C.~D. Aliprantis and O.~Burkinshaw, \emph{Positive operators}, Springer,
  Dordrecht, 2006.

\bibitem{allenby_labuschagne}
P.~D. Allenby and C.~C.~A. Labuschagne, \emph{On the uniform density of
  {$C(X)\otimes C(Y)$} in {$C(X\times Y)$}}, Indag. Math. (N.S.) \textbf{20}
  (2009), no.~1, 19--22.

\bibitem{buskes_loomis_sikorski}
G.~Buskes, B.~de~Pagter, and A.~van Rooij, \emph{The {L}oomis-{S}ikorski
  theorem revisited}, Algebra Universalis \textbf{58} (2008), no.~4, 413--426.

\bibitem{buskes_thorn_complete}
G.~Buskes and L.P. Thorn, \emph{Two results on {F}remlin's {A}rchimedean
  {R}iesz space tensor product}, arXiv:2206.06283 (July 2022).

\bibitem{fremlin_tensor_1972}
D.~H. Fremlin, \emph{Tensor products of {A}rchimedean vector lattices}, Amer.
  J. Math. \textbf{94} (1972), 777--798.

\bibitem{fremlin_measure}
\bysame, \emph{Measure theory. {V}ol. 3}, Torres Fremlin, Colchester, 2004,
  Measure algebras, Corrected second printing of the 2002 original.

\bibitem{goffman_1958}
Casper Goffman, \emph{Remarks on lattice ordered groups and vector lattices.
  {I}. {C}arath\'{e}odory functions}, Trans. Amer. Math. Soc. \textbf{88}
  (1958), 107--120.

\bibitem{grobler_2022}
J.~Grobler, \emph{Lattice tensor products in different categories of {R}iesz
  spaces}, Research Gate (July 2022).

\bibitem{jakubik_place}
J\'{a}n Jakub\'{\i}k, \emph{On vector lattices of elementary {C}arath\'{e}odory
  functions}, Czechoslovak Math. J. \textbf{55(130)} (2005), no.~1, 223--236.

\bibitem{luxemburg_riesz_1971}
W.~A.~J. Luxemburg and A.~C. Zaanen, \emph{Riesz spaces. {V}ol. {I}},
  North-Holland Publishing Co., Amsterdam-London, 1971.

\bibitem{zaanen_introduction_1997}
A.~C. Zaanen, \emph{Introduction to operator theory in {R}iesz spaces},
  Springer-Verlag, Berlin, 1997.

\end{thebibliography}

\end{document}